\documentclass[12pt]{article}
\usepackage{amsmath}

\font\tenmsb=msbm10
\font\sevenmsb=msbm7
\font\fivemsb=msbm5
\newfam\msbfam
\textfont\msbfam=\tenmsb
\scriptfont\msbfam=\sevenmsb
\scriptscriptfont\msbfam=\fivemsb

\def\mathbb#1{{\fam\msbfam\relax#1}}

\newcommand{\al}{\alpha}
\newcommand{\io}{+\infty}
\newcommand{\lm}{\left}
\newcommand{\rt}{\right}

\title{Irrationalit\'e d'au moins
un des\\ neuf nombres $\zeta(5), \zeta(7),\ldots,\zeta(21)$}
\author{Tanguy Rivoal\\ \\ Laboratoire SDAD, CNRS FRE 2271\\
D\'epartement de Math\'ematiques\\ Universit\'e de Caen, Campus II, BP 5186\\
14032 Caen C\'edex, France}
\bigskip
\date{}

\newtheorem{theo}{Th\'eor\`eme}
\newtheorem{lem}{Lemme}

\begin{document}
\maketitle
\begin{center}
\section{Introduction}
\end{center}
\vspace{0.3cm}


Le Th\'eor\`eme 2 de [BR] montre qu'il existe un entier impair $j$ tel
que $5\le j\le 169$ et $1$, $\zeta(3)$ et $\zeta(j)$ sont 
lin\'eairement ind\'ependants
sur $\mathbb{Q}\,$ : ce r\'esultat implique l'irrationalit\'e
de $\zeta(j)$ mais est bien s\^ur plus fort.
Dans cet article, nous am\'eliorons la majoration $j\le 169$ en ne
recherchant que l'irrationalit\'e de $\zeta(j)$ :
\vskip 1cm
\begin{theo} Il existe un entier impair $j$ tel que
$5\le j\le 21$ \mbox{et $\zeta(j)\not\in\,\mathbb{Q}$.}
\end{theo}
La d\'emonstration de ce th\'eor\`eme repose sur la s\'erie suivante
$$
S_{n,a}(z)=n!^{a-6}\sum_{k=1}^{\io}\frac{1}{2}\frac{\textup{d}^2}
{\textup{d}t^2}\lm\{\lm(t+\frac{n}{2}\rt)\frac{(t-n)_n^3(t+n+1)_n^3}{( 
t)_{n+1}^a}\rt\}_{|t=k}z^{-k}
$$
o\`u $z$ est un nombre complexe de module $\ge 1$ et $a$ un entier 
$\ge 6$.\\ \null

L'\'etude de $S_{n,a}(z)$, que nous \'ecrirons d\'esormais $S_n(z)$, 
est similaire \`a
celle de la s\'erie consid\'er\'ee dans [R] et [BR] :\\ \null

\noindent $\bullet$  Le Lemme 1 montre que, si $a$ est pair, la 
s\'erie $S_n(1)$
s'\'ecrit comme une combinaison lin\'eaire (\`a coefficients rationnels) de
1 et des $\zeta(j)$ pour $j$ impair, $j\in\{5,\ldots,a+2\}$.\\ \null

\noindent $\bullet$ Le Lemme 2 d\'etermine un d\'enominateur commun aux
coefficients de cette combinaison lin\'eaire.\\ \null

\noindent $\bullet$ L'estimation du comportement de
$\vert S_n(1)\vert^{1/n}$ est d\'elicate puisqu'une expres\-sion
int\'egrale de type Beukers [Be] n'est pas connue pour
$S_n(1)$. N\'ean\-moins, en suivant Nesterenko [Ne], le Lemme 4 montre
que $S_n(1)$ peut s'\'ecrire comme la partie r\'eelle d'une 
int\'egrale complexe :
le comportement asymptotique de cette int\'egrale
est alors d\'etermin\'e au Lemme 5 par la m\'ethode du col (Lemme 3).\\ \null

\noindent $\bullet$ Enfin, il n'y a pas lieu ici de borner la hauteur des
coefficients de la combinaison : cela n'est n\'ecessaire
que pour l'ind\'ependance lin\'eaire.
\vskip 1cm
\noindent \textbf{Remerciements} L'auteur tient \`a remercier F. Amoroso et D. Essouabri
pour leurs conseils qui ont permis d'am\'eliorer une pr\'ec\'edente version.

\section{R\'esultats auxiliaires}
Posons
$$
R_{n}(t)=n!^{a-6}\lm(t+\frac{n}{2}\rt)\frac{(t-n)_n^3(t+n+1)_n^3}{(t)_{n+1}^a},
$$
$D_{\lambda}=\frac{1}{\lambda!}\lm(\frac{\textup{d}}
{\textup{d}t}\rt)^{\lambda}$ et $c_{l,j,n}=D_{a-l}
(R_n(t)(t+j)^a)_{|t=-j}$ : on a
alors la d\'ecomposition en \'el\'ements simples
\begin{eqnarray}
R_n''(t)=\sum_{l=1}^a\sum_{j=0}^n
\frac{l(l-1)c_{l,j,n}}{(t+j)^{l+2}}\;.
\end{eqnarray}
D\'efinissons \'egalement les polyn\^omes \`a coefficients rationnels
\begin{eqnarray}
P_{0,n}(z)=-\sum_{l=1}^a\sum_{j=1}^n
\sum_{k=1}^{j}\frac{l(l-1) c_{l,j,n}}{2k^{l+2}} z^{j-k}
\quad\text{et}\quad P_{l,n}(z)=\sum_{j=0}^n c_{l,j,n} z^j\;.
\end{eqnarray}
o\`u $l\in\{1,\ldots,a\}$
\begin{lem}
Pour tout $z\in\mathbb{C}$, $|z|>1$, on a
$$
S_n(z)=P_{0,n}(z)+\sum_{l=1}^{a} \frac{l(l-1)}{2}
P_{l,n}(z)\textup{Li}_{l+2}(1/z)
$$
et $P_{1,n}(1)=0$. De plus, si $a$ est pair, alors pour tout $n\ge 0$ 
et pour tout entier
pair $l\in\{2,\ldots,a\}$, on a $P_{l,n}(1)=0$ et donc
$$
S_n(1)=P_{0,n}(1)+\sum_{j=2}^{a/2} j(2j-1)P_{2j-1,n}(1)\zeta(2j+1)\;.
$$
\end{lem}

\noindent \textbf{D\'emonstration}\\
\noindent De la d\'ecomposition (1) de $R_{n}(t)$, on d\'eduit
que si $|z|>1$
\begin{eqnarray*}
S_n(z) &=&
\sum_{l=1}^a\sum_{j=0}^n \frac{l(l-1) c_{l,j,n}}{2}
\sum_{k=1}^{\io}\frac{z^{-k}}{(k+j)^{l+2}}\\
&=&\sum_{l=1}^a\sum_{j=0}^n \frac{l(l-1) c_{l,j,n}}{2}
z^j\lm(\sum_{k=1}^{\io}\frac{1}{k^{l+2}}z^{-k}-\sum_{k=1}^{j}\frac{1}{ 
k^{l+2}}z^{-k}\rt)\\
&=&P_{0,n}(z)+\sum_{l=1}^a \frac{l(l-1)}{2}P_{l,n}(z)\textup{Li}_{l+2}(1/z)\;.
\end{eqnarray*}
\noindent Comme le degr\'e total de la fraction rationnelle $R_n(t)$
est $\le -2$, on a
$$
P_{1,n}(1)=\sum_{j=0}^n \text{Res}_{t=-j}(R_n(t))=0\;.
$$
On peut r\'e\'ecrire $
c_{l,j,n}=(-1)^{a-l}D_{a-l}(\Phi_{n,j}(x))_{\vert x=j}
$
o\`u
\begin{eqnarray*}
\Phi_{n,j}(x)
=n!^{a-6}\lm(\frac{n}{2}-x\rt)\frac{(-x-n)_{n}(-x+n+1)_{n}}{(-x)_{n+1} 
^a}(j-x)^a\;.
\end{eqnarray*}
On a
\begin{eqnarray}
\Phi_{n,n-j}(n-x)=
n!^{a-6}\lm(x-\frac{n}{2}\rt)
\frac{(x-2n)_{n}(x+1)_{n}}{(x-n)_{n+1}^a}
(x-j)^a\;.
\end{eqnarray}
En appliquant l'identit\'e
$(\al)_l=(-1)^l(-\al-l+1)_l$ aux trois symboles de Pochhammer de (3),
on obtient
\begin{eqnarray*}
\lefteqn{\Phi_{n,n-j}(n-x)}\\
&=&-n!^{a-6}\lm(\frac{n}{2}-x\rt)\frac{(-1)^{n}(-x+n+1)_{n}
(-1)^{n}(-x-n)_{n}}{(-1)^{(n+1)a}(-x)_{n+1}^a}(-1)^a(j-x)^a
\\
&=& (-1)^{na+1}\Phi_{n,j}(x)\;.
\end{eqnarray*}
Donc pour tout $k\geq 0$,
$$
\Phi_{n,n-j}^{(k)}(n-x)=(-1)^{k+na+1}\Phi_{n,j}^{(k)}(x)\;.
$$\\
En particulier, avec $k=a-l$ et $x=j$, on a
$$
c_{l,n-j,n}=(-1)^{a(n+1)+l+1}c_{l,j,n}\;,
$$
ce qui implique la relation
$$
P_{l,n}(1)=(-1)^{(n+1)a+l+1}P_{l,n}(1)\;.
$$
Si $(n+1)a+l$ est pair, on en d\'eduit que $P_{l,n}(1)=0$.\\

\begin{lem}
Pour tout $l\in\{1,\ldots,a\}$ on a
$$
2d_n^{a-l}P_{l,n}(z)\in\,\mathbb{Z}[z]\quad\text{ et }\quad 2d_n^{a+2}
P_{0,n}(z)\in\,\mathbb{Z}[z]
$$
o\`u $d_n=\textup{ppcm}(1,2,\ldots,n)$.
\end{lem}

\noindent \textbf{D\'emonstration}\\
On \'ecrit
$
R_n(t)(t+j)^a
=F(t)^3\times G(t)^3\times H(t)^{a-6}\times I(t)
$ o\`u $I(t)=t+n/2$ et
$$
F(t)=\frac{(t-n)_n}{(t)_{n+1}}(t+j)\,,\;
G(t)=\frac{(t+n+1)_n}{(t)_{n+1}}(t+j)\,,\;
H(t)=\frac{n!}{(t)_{n+1}}(t+j)\;.
$$
D\'ecomposons $F(t)$, $G(t)$  et $H(t)$ en fractions
partielles :
$$
F(t)=1+\sum_{\substack{p=0\\p\not=j}}^n
\frac{j-p}{t+p}f_{p}\,,\quad
G(t)=1+\sum_{\substack{p=0\\p\not=j}}^n
\frac{j-p}{t+p}g_{p}\,,\quad H(t)=
\sum_{\substack{p=0\\p\not=j}}^n \frac{j-p}{t+p}h_{p}
$$
o\`u
$$
f_{p}=\frac{(-p-n)_n}{\displaystyle\prod_{\substack{h=0\\h\not=p}}^n
(-p+h)} =\frac{(-1)^n (p+1)_n}{(-1)^p p!(n-p)!}
=(-1)^{n-p}\binom{n+p}{n}\binom{n}{p}\in\mathbb{Z}\;,
$$
$$
g_{p}=\frac{(-p+n+1)_n}{\displaystyle\prod_{\substack{h=0\\h\not=p}}^n
(-p+h)}=\frac{(-1)^p(2n-p)!}{(n-p)!p!(n-p)!}=
(-1)^p\binom{2n-p}{n}\binom{n}{p} \in\mathbb{Z}
$$
et
$$
h_p=\frac{n!}{\displaystyle\prod_{\substack{h=0\\h\not=p}}^n
(-p+h)}=\frac{(-1)^p n!}{p!(n-p)!} =(-1)^p \binom{n}{p} \in\mathbb{Z}\;.
$$
On a alors pour tout entier $\lambda\geq 0$ :
$$
(D_{\lambda}F(t))_{\vert t=-j}
=\delta_{0,\lambda}+\sum_{\substack{p=0\\p\not=j}}^n
(-1)^{\lambda}\frac{j-p}{(p-j)^{\lambda+1}}f_{p}\;,
$$
$$
(D_{\lambda}G(t))_{\vert t=-j}=
\delta_{0,\lambda}+\sum_{\substack{p=0\\p\not=j}}^n
(-1)^{\lambda}\frac{j-p}{(p-j)^{\lambda+1}}g_{p}\;,
$$
$$
(D_{\lambda}H(t))_{\vert t=-j}=
\sum_{\substack{p=0\\p\not=j}}^n
(-1)^{\lambda}\frac{j-p}{(p-j)^{\lambda+1}}h_p
$$
avec $\delta_{0,\lambda}=1$ si $\lambda=0$, $\delta_{0,\lambda}=0$
si $\lambda>0$.  On a donc montr\'e que
$$
d_n^{\lambda}(D_{\lambda}F)_{\vert t=-j}\;,
\quad d_n^{\lambda}(D_{\lambda}G)_{\vert t=-j}
\quad\text{et}\quad d_n^{\lambda}(D_{\lambda}H)_{\vert t=-j}
$$
sont des entiers pour tout $\lambda \in\mathbb{N}\,$. De plus,
$2(D_{\lambda}I)_{\vert t=-j}\in\mathbb{Z}$. Gr\^ace \`a
la formule de Leibniz
\begin{multline*}
D_{a-l}(R(t)(t+j)^a)=
\sum_{\mu}(D_{\mu_1}F)(D_{\mu_2}F)(D_{\mu_3}F)\\
\quad\times(D_{\mu_4}G)(D_{\mu_5}G)(D_{\mu_6}G)
(D_{\mu_{7}}H)\cdots (D_{\mu_{a}}H)(D_{\mu_{a+1}}I)
\end{multline*}
(o\`u la somme est sur les multi-indices $\mu\in\mathbb{N}\,^{a+1}$
tels que $\mu_1+\cdots+\mu_{a+1}=a-l$), on en d\'eduit alors que
$2d_n^{a-l}c_{l,j,n}\in\mathbb{Z}\,$. Les expressions (2) des polyn\^omes
$P_{0,n}(z)$ et $P_{l,n}(z)$ permettent de conclure.


\section{D\'emonstration du Th\'eor\`eme 1}

Pour estimer $S_n(1)$, nous suivons
la d\'emarche utilis\'ee par [Ne] et [HP] qui consiste
\`a exprimer $S_n(1)$ \`a l'aide d'une
int\'egrale complexe \`a laquelle on peut appliquer la m\'ethode du col,
m\'ethode dont nous rappelons tout d'abord le principe
(voir par exemple [Co], pp. 91-94 ou [Di], pp. 279-285]).\\ \null

Soit $w$ une fonction analytique au voisinage d'un point $z_0$. On appelle
chemin de descente de $\text{Re}(w)$ en $z_0$ tout chemin
du plan issu de $z_0$ et le long duquel
$\text{Re}(w(z))$ est strictement d\'ecroissante quand $z$ s'\'eloigne de
$z_0$. Les chemins de plus grande descente de $\text{Re}(w)$ en $z_0$ sont
les chemins
tels que $\text{Re}(w)$ a (localement) la d\'ecroissance la plus rapide parmi
tous les chemins de descente : il est en fait \'equivalent
de demander que $\text{Im}(w)$ soit constante le long de
ces chemins, c'est \`a dire que la phase de $e^w$ soit stationaire.\\ \null

Supposons $w$ telle que
$w^{\prime}(z_0)=0$ et
$w^{\prime\prime}(z_0)=\vert w^{\prime\prime}(z_0) \vert e^{i\al_0}\not=0$.
Notons $\theta$ la direction d'une droite $\Delta$
passant par $z_0$, c'est-\`a-dire
$\theta=\text{arg}(z-z_0)$ o\`u $z\in\,\Delta$. Il existe exactement deux
chemins de
plus grande descente de $\text{Re}(w)$
en $z_0$, dont les directions des tangentes en $z_0$ sont
$\theta_+=\frac{\pi}{2}-\frac{\al_0}{2}$ et
$\theta_-=-\frac{\pi}{2}-\frac{\al_0}{2}$ : ces directions
\textit{critiques}
sont oppos\'ees. Il peut s'av\'erer difficile de d\'eterminer exactement
les chemins
de plus grande descente. On peut s'affranchir de ce probl\`eme
en consid\'erant n'importe quelle direction $\theta$ en $z_0$
telle que $\cos(\al_0+2\theta)<0$ : au voisinage de $z_0$,
$$
w(z)=w(z_0)+\frac{1}{2} w^{\prime\prime}(z_0) (z-z_0)^2+O((z-z_0)^3)
$$
et sur un chemin $L$ dont les deux directions en $z_0$ v\'erifient la
condition ci-dessus,
on a alors $\text{Re}(\frac{1}{2} w^{\prime\prime}(z_0)(z-z_0)^2)<0$
et $\text{Re}(w)$ admet un maximum local en $z_0$ le long de $L$.
Convenons de dire qu'un chemin $L$ est \textit{admissible} en $z_0$
si les deux directions
$\theta$ en $z_0$ v\'erifient $\cos(\al_0+2\theta)<0$
et si $\text{Re}(w(z_0))$ est le maximum \textit{global} de $\text{Re}(w)$
le long de $L$.\\

\begin{lem}[M\'ethode du col]
Soit $g$ et $w$ deux fonctions analytiques dans un ouvert simplement
connexe $\mathcal{D}$
du plan. Supposons qu'il existe $z_0\in\,\mathcal{D}$ tel que
$w^{\prime}(z_0)=0$ et
$w^{\prime\prime}(z_0)= \vert w^{\prime\prime}(z_0)\vert
e^{i\al_0}\not=0$. Si $L$ est un chemin inclus dans $\mathcal{D}$ et
admissible en $z_0$, alors
\begin{eqnarray}
\int_L g(z)e^{nw(z)}\textup{d}z \sim g(z_0)
\sqrt{\frac{2\pi}{n\vert w^{\prime\prime}(z_0) \vert}}
\,e^{i(\pm\frac{\pi}{2}-\frac{\al_0}{2})}e^{nw(z_0)}\quad (n\to\io)
\end{eqnarray}
o\`u le choix de $\pm$ d\'epend de l'orientation de $L$.
De plus, cette estimation est encore valable si $L$ est un chemin que l'on peut
d\'eformer en un chemin admissible en $z_0$.
\end{lem}

Nous appliquons maintenant cette m\'ethode \`a
l'estimation asymptotique de $S_n(1)$. Consid\'erons l'int\'egrale complexe
$$
J_n(u)=\frac{n}{2i\pi}\int_L R_n(nz)\lm(\frac{\pi}{\sin(n\pi z)}\rt)^3
e^{nuz}\textup{d}z
$$
o\`u $u$ est un nombre
complexe tel que $\text{Re}(u)\le 0$ et $\vert \text{Im}(u) \vert\le 3\pi$,
$L$ est une droite verticale orient\'ee de $+i\infty$ \`a $-i\infty$ 
et contenue
dans la bande $0<\text{Re}(z)<1$, ce qui assure que
l'int\'egrale $J_n(u)$ converge.\\

\begin{lem}
Dans ces conditions, on a\\ \null
\textit{i)}
\begin{multline*}
J_n(u)=\frac{(-1)^n n^2}{2i\pi} n!^{a-6}\\
\times\int_L\lm(z+\frac{1}{2}\rt)
\frac{\Gamma(nz)^{a+3}\Gamma(n-nz+1)^3
\Gamma(nz+2n+1)^3}{\Gamma(nz+n+1)^{a+3}}e^{nuz}\textup{d}z\;.
\end{multline*}
\textit{ii)}
$$
S_n(1)=\textup{Re}\,(J_n(i\pi))\;.
$$
\end{lem}

\noindent \textbf{D\'emonstration}\\
\noindent \textit{i)} Comme $(\al)_n=\Gamma(\al+n)/\Gamma(\al)$ et
$(t-n)_n^3=(-1)^{n}(1-t)_n^3$, on a
\begin{eqnarray*}
R_n(t)&=&(-1)^n n!^{a-6}\lm(t+\frac{n}{2}\rt)\frac{(1-t)_n^3 
(t+n+1)_n^3}{(t)_{n+1}^a}\\
      &=&(-1)^n 
n!^{a-6}\lm(t+\frac{n}{2}\rt)\frac{\Gamma(n-t+1)^3\Gamma(t+2n+1)^3
\Gamma(t)^a}{\Gamma(1-t)^3\Gamma(t+n+1)^{3}\Gamma(t+n+1)^{a}}\;.
\end{eqnarray*}
De plus, la formule des compl\'ements $\Gamma(t)\Gamma(1-t)=\pi/\sin(\pi t)$
(pour $t\not\in\,\mathbb{Z}$) implique que
$$
R_n(t)\lm(\frac{\pi}{\sin \pi t}\rt)^3=
(-1)^n n!^{a-6}\lm(t+\frac{n}{2}\rt)\frac{\Gamma(t)^{a+3}\Gamma(n-t+1)^3
\Gamma(t+2n+1)^3}{\Gamma(t+n+1)^{a+3}}\;.
$$
On a donc
\begin{multline*}
\int_{L^{\prime}} R_n(t)\lm(\frac{\pi}{\sin(\pi t)}\rt)^3 e^{ut}\textup{d}t\\=
(-1)^n n!^{a-6}\int_{L^{\prime}}\lm(t+\frac{n}{2}\rt)\frac{\Gamma(t)^{a+3}
\Gamma(n-t+1)^3\Gamma(t+2n+1)^3}{\Gamma(t+n+1)^{a+3}}e^{ut}\textup{d}t
\end{multline*}
o\`u $L^{\prime}$ est une droite verticale quelconque contenue dans
$0<\text{Re}(t)<n$.
Le changement de variable $t=nz$ et le th\'eor\`eme de Cauchy justifient que
\begin{multline*}
J_n(u)=\frac{(-1)^n n^2}{2i\pi} n!^{a-6}\\
\times\int_L\lm(z+\frac{1}{2}\rt)
\frac{\Gamma(nz)^{a+3}\Gamma(n-nz+1)^3
\Gamma(nz+2n+1)^3}{\Gamma(nz+n+1)^{a+3}}e^{nuz}\textup{d}z\;.
\end{multline*}
\noindent \textit{ii)} Soit $c\in]0,n[$ et soit 
$T\in\frac{1}{2}+\mathbb{Z}$ tel que
$T>n+1$.  Consid\'erons le contour rectangulaire $\mathcal{R}_T$ 
orient\'e dans le
sens direct, de sommets $c \pm i T$ et $T\pm i T$ : la
fonction $F(t,u)=R_n(t)(\pi/\sin(\pi t))^3 e^{ut}$ est m\'eromorphe dans le
demi-plan
$\text{Re}(t)>0$ et ses p\^oles sont les entiers $k\ge n+1$. En appliquant le
th\'eor\`eme des r\'esidus, il d\'ecoule que
$$
\frac{1}{2i\pi}\int_{\mathcal{R}_T} F(t,u) \textup{d}t=\sum_{k=n+1}^{[T]}
\text{Res}_{t=k}(F(t,u))\;.
$$
o\`u
$$
\text{Res}_{t=k}(F(t,u))=\frac{\pi^2+u^2}{2}R_n(k)(-e^u)^k
+uR_n'(k)(-e^u)^k+\frac{1}{2}R_n''(-e^u)^k\;.
$$
Sur les trois c\^ot\'es $[c-iT,T-iT]$, $[T-iT,T+iT]$ et $[T+iT,c+iT]$,
on a $R_n(t)=O(T^{-2})$.\\ \null

\noindent Sur $[T-iT,T+iT]$, en posant $t=T+iy$, on a
$$
\sin(\pi t)=(-1)^N \cosh(\pi y)
$$
et donc $|\sin(\pi t)|\ge \frac{1}{2}e^{\pi |y|}$. Comme
$\lm|e^{ut}\rt|=e^{\text{Re}(u)T-\text{Im}(u)y}$, on en d\'eduit que
$$
R_n(t)\lm(\frac{\pi}{\sin(\pi t)}\rt)^3 e^{ut}
=O\lm(T^{-2}e^{\text{Re}(u)T}e^{-(\text{Im}(u)y+3\pi|y|)}\rt)
=O\lm(T^{-2}\rt)
$$
puisque $\text{Re}(u)\le 0$ et $|\text{Im}(u)|\le 3\pi$.\\ \null

\noindent De fa\c{c}on similaire, sur les deux c\^ot\'es $[c-iT,T-iT]$ et
$[T+iT,c+iT]$,
en posant $t=x\pm iT$ avec $x>0$, on a
$$
2 i \sin(\pi t)=e^{\mp\pi T}e^{i\pi x}-e^{\pm\pi T}e^{-i\pi x}
$$
et donc $|\sin(\pi t)|\ge |\sinh(\pi T)|\gg e^{\pi T}$. Comme
$\lm|e^{ut}\rt|=e^{\text{Re}(u)x-\text{Im}(u)T}$, on en d\'eduit que

$$
R_n(t)\lm(\frac{\pi}{\sin(\pi t)}\rt)^3 e^{ut}
=O\lm(T^{-2}e^{\text{Re}(u)x}e^{-(\text{Im}(u)T+3\pi T)}\rt)=
O\lm(T^{-2}\rt)\;.
$$
Donc
\begin{eqnarray*}
J_n(u)&=&\frac{1}{2i\pi}
\int_{c+i\infty}^{c-i\infty} F(t,u) \textup{d}t=\lim_{T\to\io} \frac{1}{2i\pi}
\int_{\mathcal{R}_T} F(t,u)\textup{d}t\\
&=&\sum_{k=n+1}^{\io} \text{Res}_{t=k}(F(t,u))\\
&=&\sum_{k=n+1}^{\io} \lm(\frac{\pi^2+u^2}{2}R_n(k)(-e^u)^k
+uR_n'(k)(-e^u)^k+\frac{1}{2}R_n''(k)(-e^u)^k \rt)\,.
\end{eqnarray*}
En particulier,
$$
J_n(i\pi)=\sum_{k=n+1}^{\io} \lm(i\pi R_n'(k)+\frac{1}{2}R_n''(k)\rt)
$$
et donc $S_n(1)=\text{Re}(J_n(i\pi))$.\\ \null

Nous utilisons maintenant la formule de Stirling sous la forme suivante
$$
\Gamma(z)=\sqrt{\frac{2\pi}{z}}\lm(\frac{z}{e}\rt)^z\lm(1+O\lm(\frac{1 
}{|z|}\rt)\rt)
$$
o\`u $\vert z\vert \to\infty$, $\vert \text{arg}(z) \vert < \pi$ et o\`u
les fonctions
$\sqrt{z}$ et $z^z=e^{z\log(z)}$ sont d\'efinies avec la d\'etermination
principale du logarithme. Sur la droite $L$, les quantit\'es
$|nz|$, $|n-nz+1|$, $|nz+2n+1|$ et $|nz+n+1|$ sont
\'equivalentes \`a des multiples constants de $n$,
d'o\`u
\begin{eqnarray}
J_n(i\pi)=i(-1)^{n+1}(2\pi)^{\frac{a}{2}-1}n^{\frac{a}{2}-4}\int_L 
g(z)e^{nw(z)}
\lm(1+O\lm(\frac{1}{n}\rt)\rt)
\textup{d}z
\end{eqnarray}
avec
$$
g(z)=\frac{\sqrt{z+1}^{\,a+3}}
{\sqrt{z}^{\,a+3}\sqrt{1-z}^{\,3}\sqrt{z+2}^{\,3}}
$$
et
\begin{eqnarray*}
w(z)&=&(a+3)z\log(z)-(a+3)(z+1)\log(z+1)\\
&&\qquad +3(1-z)\log(1-z)+3(z+2)\log(z+2)+i\pi z\;,
\end{eqnarray*}
les diff\'erentes fonctions racines et logarithmes
de $g$ et $w$ \'etant de nouveau
d\'efinies \`a l'aide de la d\'etermination
principale du logarithme. L'expression (5)
de $J_n(i\pi)$ se pr\^ete maintenant
\`a une estimation par la m\'ethode du col.\\ \null

Dor\'enavant, nous supposons $a=20$. Alors
$$
w^{\prime}(z)=23\log(z)-23\log(z+1)+3\log(z+2)-3\log(1-z)+i\pi
$$
et l'\'equation $w^{\prime}(z)=0$ poss\`ede une seule solution $z_0$ 
v\'erifiant
$0<\text{Re}(z_0)<1$ :
$$
z_0=x_0+i\,y_0\approx 0,9922341203-i\,0,01200539829\;.
$$
On a
$$
w(z_0)\approx -22,02001640+i\,3,104408624$$
et
$$
w^{\prime\prime}(z_0)\approx 216,7641546e^{-i\,0.9471277165}\;.
$$
On constate que $\theta=\pi/2$ et $\theta=-\pi/2$ v\'erifient
$\cos(\al_0+2\theta)<0$. Montrons que la droite $L:\,\text{Re}(z)=x_0$
est admissible, c'est \`a dire que $\text{Re}(w)$ admet un maximum global
en $z_0$ le long de $L$. Posons $f(y)=\frac{\partial 
\text{Re}(w)}{\partial y}(x_0+iy)$~; donc
\begin{eqnarray*}
f(y)&=&-\text{Im}(w^{\prime})(x_0+iy)\\
&=&-23\,\text{arg}(x_0+iy)+23\,\text{arg}(x_0+1+iy)\\
&&\;-3\,\text{arg}(x_0+2+iy)+3\,\text{arg}(1-x_0-iy)-\pi\;.
\end{eqnarray*}
On a
$$
\lim_{y\to -\infty}f(y)=2\pi\quad\hbox{et}\quad\lim_{y\to +\infty}f(y)=-4\pi\;.
$$
Par ailleurs, $\text{arg}(z)=\arctan\lm(\frac{\text{Im}(z)}{\text{Re}(z)}\rt)$
pour $\text{Re}(z)>0$, d'o\`u
\begin{eqnarray*}
\frac{\textup{d}f}{\textup{d}y}&=&-\frac{23x_0}{x_0^2+y^2}+\frac{23(x_ 
0+1)}{(x_0+1)^2+y^2}
-\frac{3(x_0+2)}{(x_0+2)^2+y^2}-\frac{3(1-x_0)}{(1-x_0)^2+y^2}\\
&=&\frac{N(y^2)}{(x_0^2+y^2)((x_0+1)^2+y^2)((x_0+2)^2+y^2)((1-x_0)^2+y^2)}\;,
\end{eqnarray*}
o\`u l'on a not\'e
\begin{eqnarray*}
N(t)\hskip -7pt&=&\hskip -7pt14t^3+2(7x_0^2+7x_0+44)t^2
+2(-7x_0^4-14x_0^3-124x_0^2-117x_0+37)t
\\&&+2(-7x_0^5-21x_0^4+16x_0^3+67x_0^2-9)x_0\;.
\end{eqnarray*}
On v\'erifie que $N(t)$ a une seule racine dans $[\,0,\io\,[$. Donc
$f(y)$ ne s'annule que pour $y=y_0$. La fonction
$y\to \text{Re}(w(x_0+iy))$ est donc strictement
croissante sur $]-\infty, y_0\,]$, puis strictement
d\'ecroissante sur $[\,y_0,+\infty[$. En cons\'equence,
la droite $L:\,\text{Re}(z)=x_0$ est admissible en $z_0$ pour $\text{Re}(w)$.\\

\begin{lem}
On a :
$$
J_n(i\pi)\sim c_0(-1)^{n+1}n^{11/2}e^{nw(z_0)}\quad (n\to\io)
$$
o\`u 
$c_0=g(z_0)(2\pi)^{19}\sqrt{2\pi/|w^{\prime\prime}(z_0)|}e^{-i\al_0/2} 
\not=0$. De plus, il
existe une suite d'entiers $\varphi(n)$ telle que
$$
\limsup_{n\to\io}|S_{\varphi(n)}(1)|^{1/\varphi(n)}=e^{\text{Re}(w(z_0))}
$$

\end{lem}
\noindent \textbf{D\'emonstration}\\
L'estimation de $J_n(i\pi)$ r\'esulte de l'estimation g\'en\'erale (4),
appliqu\'ee \`a (5) et \`a la droite admissible 
\mbox{$L:\,\text{Re}(z)=x_0$.} Pour montrer la
derni\`ere affirmation, notons $c_0=r\,e^{i\beta}$ et 
$v_0=\text{Im}(w(z_0))$, de sorte que
\begin{eqnarray*}
S_n(1)\hskip -8pt&=&\hskip -8pt\text{Re}(J_n(i\pi))\\
\rule{0mm}{6mm}\hskip -8pt&=&\hskip -8pt 
r(-1)^{n+1}n^{11/2}e^{n\text{Re}(w(z_0))}
(\text{Re}(u_n)\cos(n v_0+\beta)-\text{Im}(u_n)\sin(n v_0+\beta))
\end{eqnarray*}
o\`u $u_n$ est une suite de nombres complexes qui converge vers $1$. 
Remarquons que $v_0\approx 3,104$ n'est
pas un multiple entier de $\pi$ et donc il existe une suite d'entiers 
$\varphi(n)$ telle que
$\cos(\varphi(n)v_0+\beta)$ converge vers une limite $l\neq0$. On en 
d\'eduit que
$$
\lim_{n\to\io}(\text{Re}(u_{\varphi(n)})\cos(\varphi(n)v_0+\beta)
-\text{Im}(u_{\varphi(n)})\sin(\varphi(n)v_0+\beta))=l\neq0
$$
et donc
$$
\lim_{n\to\io}|S_{\varphi(n)}(1)|^{1/\varphi(n)}=e^{\text{Re}(w(z_0))}\;.
$$
\vskip 0.5cm
\noindent \textbf{D\'emonstration du Th\'eor\`eme 1}\\
Posons $p_{0,n}=2d_{n}^{22}P_{0,n}(1)$ et 
$p_{l,n}=2l(2l-1)d_{n}^{22}P_{2l-1,n}(1)$
pour $l\in\{2,\dots,10\}$ : le Lemme 2 implique que ce sont des entiers.
D\'efinissons \'egalement $\ell_n=2d_{n}^{22}S_{n}(1)$ : le Lemme 1 montre que
$$
\ell_n=p_{0,n}+\sum_{l=2}^{10} p_{l,n}\zeta(2l+1)\;.
$$
Enfin, d'apr\`es le Th\'eor\`eme des nombres premiers, $d_n=e^{n+o(n)}$.
Le Lemme 5 montre que
$$
\lim_{n\to\io}|\ell_{\varphi(n)}|^{1/{\varphi(n)}}\approx 
e^{-0,02}\in\,]\,0,1\,[\;,
$$
ce qui prouve le Th\'eor\`eme 1.\\ \null

\vspace{0.2cm}
\centerline{\textbf{R\'ef\'erences}}
\vspace{0.2cm}

\noindent [BR] K. Ball et T. Rivoal, \textit{Irrationalit\'e d'une infinit\'e
de valeurs de la fonction z\^eta aux entiers impairs}, soumis.\\ \null

\noindent [Be] F. Beukers, \textit{A note on the irrationality of
$\zeta(2)$ and $\zeta(3)$},  Bull. Lond. Math. Soc.
\textbf{11}, no. 33, 268-272 (1978).\\ \null

\noindent [Co] E. T. Copson, \textit{Asymptotic expansions},
Cambridge University Press (1967).\\ \null

\noindent [Di] J. Dieudonn\'e, \textit{Calcul infinit\'esimal},
Collection "M\'ethodes", Hermann (1980).\\ \null

\noindent [HP] T. G. Hessami Pilerhood, \textit{Linear independence of vectors
with\break polylogarithmic coordinates}, Mosc. Univ. Math. Bull. \textbf{54},
no. 6, 40-42 (1999).\\ \null

\noindent [Ne] Yu.V. Nesterenko, \textit{A few remarks on $\zeta(3)$},
Math. Notes, \textbf{59}, no. 6,  625-636 (1996).\\ \null

\noindent [R] T. Rivoal, \textit{La fonction Z\^eta de Riemann prend une
infinit\' e de valeurs irrationnelles aux entiers impairs}, C. R. Acad. Sci.
Paris \textbf{331}, 267-270 (2000).\\ \null

\end{document}